\documentclass[a4paper]{amsart}

\usepackage{amsmath}
\usepackage{amssymb}
\usepackage{amsthm}
\usepackage{graphicx}

\def\R{\ensuremath{\mathbb{R}}}
\def\Z{\ensuremath{\mathbb{Z}}}
\def\N{\ensuremath{\mathbb{N}}}
\def\a{\ensuremath{\boldsymbol{a}}}
\def\b{\ensuremath{\boldsymbol{b}}}
\def\c{\ensuremath{\boldsymbol{c}}}
\def\seqN{\ensuremath{\N^\N}}

\newcommand{\floor}[1]{{\ensuremath{\lfloor#1\rfloor}}}
\newcommand{\ceil}[1]{{\ensuremath{\lceil#1\rceil}}}
\def\epsilon{\varepsilon}
\def\phi{\varphi}

\newtheorem{thm}{Theorem}
\newtheorem{pro}[thm]{Proposition}
\newtheorem{lem}[thm]{Lemma}
\newtheorem{cor}[thm]{Corollary}
\theoremstyle{definition}
\newtheorem{defn}{Definition}
\theoremstyle{remark}
\newtheorem{rem}{Remark}
\newtheorem{exmp}{Example}
\theoremstyle{plain}

\def\BExe{\begin{exmp}}
\def\EExe{\end{exmp}}
\def\BRem{\begin{rem}}
\def\ERem{\end{rem}}
\def\BThe{\begin{thm}}
\def\EThe{\end{thm}}
\def\BDef{\begin{defn}}
\def\EDef{\end{defn}}
\def\BPro{\begin{pro}}
\def\EPro{\end{pro}}
\def\BLem{\begin{lem}}
\def\ELem{\end{lem}}
\def\BCor{\begin{cor}}
\def\ECor{\end{cor}}
\def\BProof{\begin{proof}}
\def\EProof{\end{proof}}

\begin{document}

\author{Samuel Nicolay}
\address{Corresponding author. Email: S.Nicolay@ulg.ac.be. Phone: +32(0)43669433. Fax: +32(0)43669547.}
\author{Laurent Simons}
\address{Universit\'e de Li\`ege, Institut de Math\'ematique, Grande Traverse, 12, B\^atiment B37, B-4000 Li\`ege (Sart-Tilman), Belgium.}
\title{About the multifractal nature of Cantor's bijection}
\date{\today}

\begin{abstract}
In this note, we investigate the regularity of Cantor's one-to-one mapping between the irrational numbers of the unit interval and the irrational numbers of the unit square. In particular, we explore the fractal nature of this map by showing that its H\"older regularity lies between 0.35 and 0.72 almost everywhere (with respect to the Lebesgue measure).
\end{abstract}
\maketitle 

{\bf MSC}{: 26A30, 11K50.}\hfil

\section{Introduction}

In 1878 \cite{cantor}, Cantor proved that there exists a one-to-one correspondence between the points of the unit line segment $[0,1]$ and the points of the unit square $[0,1]^2$ (repeated application of this result gives a bijective correspondence between $[0,1]$ and $[0,1]^n$, where $n$ is a natural number). About this discovery he wrote to Dedekind: ``Je le vois, mais je ne le crois pas~!'' (``I see it, but I don't believe it!'') \cite{wallace,gouvea}. Since this application is defined via continued fractions, it is very hard to have any intuition about its regularity. When looking at its definition or at the graphical representation of each component (given for the first time here), it is not hard to convince oneself that the behavior of such a function is necessarily ``erratic''; however, its (H\"older-)regularity has never been considered.

The set of the natural numbers is denoted by $\N$ (and does not contain $0$). We set $E=[0,1]$, denote by $D$ the rational numbers of $E$ and set $I = E\setminus D$. The set of the (infinite) sequences of natural numbers is denoted $\seqN$; since this space is a countable product of metric spaces, if $\a=(a_j)_{j\in\N}$ and $\b=(b_j)_{j\in\N}$ are two elements of $\seqN$, we define the usual distance
\[
 d(\a,\b)= \sum_{j=1}^\infty 2^{-j} \frac{|a_j-b_j|}{|a_j-b_j|+1}.
\]
We will implicitly consider that $\seqN$ is equipped with this distance, while $E$, $D$ and $I$ are endowed with the Euclidean distance.
\BRem
Considering $\a$ and $\b$ as two infinite words on the alphabet $\N$ \cite{lothaire}, one can also use the following ultrametric distance on $\seqN$: if $\a =(a_j)_{j\in\N}$ and $\b=(b_j)_{j\in\N}$ both belong to $\seqN$, let $\a\wedge \b$ denote the longest common prefix of $\a$ and $\b$, so that the length $|\a\wedge \b|$ of this prefix is equal to the lowest natural number $j$ such that $a_j\not=b_j$ minus $1$. A distance between $\a$ and $\b$ is given by
\[
 d'(\a,\b)=\left\{\begin{tabular}{ll}
                $0$               & if $\a=\b$ \\
                $2^{-|\a\wedge \b|}$ & if $\a\not=\b$
               \end{tabular}\right..
\]
The following relations hold:
\[
 2^{-2} d\le d'<d.
\]
\ERem
For the sake of completeness, let us recall the following result.
\BPro
The space $\seqN$ (endowed with the distance defined above) is a separable complete metric space.
\EPro
\BProof
If $\seqN_n$ denotes the set $\{\a=(a_j)_{j\in\N} \in \seqN: a_j=1\,\forall j>n\}$ ($n \in \N$), one directly checks that $\cup_{n\in\N} \seqN_n$ is dense in $\seqN$. Moreover, if $\a_j$ is a Cauchy sequence of $\seqN$, there exists a subsequence $\b_j$ such that $d(\b_j,\b_{j+1})<2^{-j}$ for any $j\in\N$. One easily checks that $\b_j$ converges to $\a_0 \in\seqN$ as $j$ tends to infinity, where $a_{0,k}=b_{k,k}$ ($k\in\N$). 
\EProof

In this note, we first recall the construction of Cantor's bijection between $I$ and $I^2$ based on continued fractions and give, as far as we know for the first time, a graphical representation of the two components of this map. We then construct an homeomorphism between $I$ and $\seqN$ to show that Cantor's bijection between $I$ and $I^2$ is an homeomorphism and that any extension of this mapping to $E$ is necessarily discontinuous at every rational number. We also investigate the multifractal nature of this function. It is well known that most of the ``historical'' space filling functions are monoH\"older with H\"older exponent equal to $1/2$ \cite{jaffard-nicolay,jaffard-nicolay2}; here we show that for Cantor's bijection, almost every point of $I$ (with respect to the Lebesgue measure) is associated to an H\"older exponent which belongs to an interval containing $1/2$ (more precisely, this interval is bounded by $0.35$ and $0.72$). All the obtained results strongly rely on the theory of the continued fractions (see e.g.\ \cite{khintchine}).

\section{Definitions}

\subsection{Continued fractions}
Let us first recall the basic facts about the continued fractions \cite{khintchine}. Here, we state the results for $E$, but they can be easily extended to the whole real line.

Let $\a=(a_j)_{j\in\{1,\ldots,n\}}$ be a finite sequence of positive real numbers ($n\in \N$); the expression $[a_1,\ldots,a_n]$ is recursively defined as follows:
\[
 [a_1]=1/a_1
 \quad\text{and}\quad
 [a_1,\ldots,a_m]=\frac{1}{a_1 +[a_2,\ldots,a_m]},
\]
for any $m\in\{2,\ldots,n\}$. If $\a\in \N^n$, we say that $[a_1,\ldots,a_n]$ is a (simple) finite continued fraction.
\BPro
For any $\a\in\N^n$ ($n\in\N$), $[a_1,\ldots,a_n]$ belongs to $D$. Conversely, for any $x\in D$, there exists a natural number $n$ and a sequence $\a\in\N^n$ such that $x=[a_1,\ldots,a_n]$.
\EPro
The representation of a rational number as a continued fraction is not unique, as shown by the following remark; this will be used in the proof of Proposition~\ref{pro:extension}.
\BRem
 If $\a\in\N^n$ ($n\in\N$) is such that $a_n>1$, one has
\[
 [a_1,\ldots,a_n]=[a_1,\ldots,a_n-1,1].
\]
\ERem

Let us now define the notion of convergent.
 For $\a\in\N^n$ ($n\in\N$) and each integer $j\in\{-1,\ldots,n\}$, let us define $p_j(\a)$ and $q_j(\a)$ by setting $p_{-1}(\a)=1$, $q_{-1}(\a)=0$, $p_0(\a)=0$, $q_0(\a)=1$ and, for $j\in\N$,
\[
 \left\{\begin{array}{c}
  p_j(\a)=a_j p_{j-1}(\a) +p_{j-2}(\a) \\
  q_j(\a)=a_j q_{j-1}(\a) +q_{j-2}(\a)
 \end{array}\right..
\]
The quotient $p_j(\a)/q_j(\a)$ is called the convergent of order $j$ of $\a$. They are intimately related to the continued fractions.
\BPro\label{pro:convergents}
 Let $[a_1,\ldots,a_n]$ ($n\in\N$) be a continued fraction; the section $[a_1,\ldots,a_j]$, with $j\le n$, is equal to $p_j(\a)/q_j(\a)$. Furthermore, we have, for any $j\ge 1$,
\[
 q_j(\a) p_{j-1}(\a) - p_j(\a)q_{j-1}(\a)= (-1)^j,
\]
and, for any $j\ge 2$,
\[
 q_j(\a) p_{j-2}(\a) - p_j(\a) q_{j-2}(\a)= (-1)^{j-1} a_j.
\]
As a consequence, one has, for any $j\ge 2$,
\[
 \frac{p_{j-1}(\a)}{q_{j-1}(\a)}- \frac{p_j(\a)}{q_j(\a)}= \frac{(-1)^j}{q_j(\a)q_{j-1}(\a)},
\]
and, for any $j\ge 3$,
\[
 \frac{p_{j-2}(\a)}{q_{j-2}(\a)}- \frac{p_j(\a)}{q_j(\a)}= \frac{(-1)^{j-1} a_j}{q_j(\a)q_{j-2}(\a)}.
\]
\EPro

Of course, one can define the numbers $p_j(\a)$ and $q_j(\a)$ for an element $\a$ of $\seqN$. The convergents allow to introduce the notion of infinite continued fraction, thanks to the following trivial result.
\BCor
 For any $\a\in\seqN$, the sequences
\[
 (\frac{p_{2j-1}(\a)}{q_{2j-1}(\a)})_{j\in\N}
 \quad\text{and}\quad
 (\frac{p_{2j}(\a)}{q_{2j}(\a)})_{j\in\N}
\]
are two adjacent sequences, $p_{2j}(\a)/q_{2j}(\a)$ being the increasing one.
\ECor
This shows that for any $\a \in \seqN$, the sequence $x_j=[a_1,\ldots,a_j]$ converges. The limit is called an infinite continued fraction and is denoted $[a_1,\ldots]$. If the real number $x\in E$ is equal to $[a_1,\ldots]$, we say that $[a_1,\ldots]$ is a continued fraction corresponding to $x$.
The following result states that the continued fraction is an instrument for representing  the real numbers (of $E$).
\BThe\label{thm:expr cf}
We have $x\in I$ if and only if there exists an infinite continued fraction corresponding to $x$; moreover, this infinite continued fraction is unique.
\EThe
\BPro\label{pro:quotient avec reste}
If $x\in E$ can be written as $x=[a_1,\ldots,a_n,r_{n+1}]$, with $a_1,\ldots,a_n \in\N$ and $r_{n+1}\in [1,\infty)$, the following relation holds:
\[
 x=\frac{p_n(\a) r_{n+1} +p_{n-1}(\a)}{q_n(\a)r_{n+1}+q_{n-1}(\a)},
\]
with $\a=(a_j)_{j\in\{1,\ldots,n\}}$.
\EPro

A sequence $\a \in \seqN$ is ultimately periodic of period $k\in\N$ if there exists $J$ such that $a_{j+k}=a_j$ for any $j\ge J$. In this case, the corresponding continued fraction $[a_1,\ldots]$ is also called ultimately periodic of period $k$.
The quadratic numbers (of $E$) are characterized by their corresponding continued fractions.
\BThe
An element of $I$ is a quadratic number if and only if the corresponding continued fraction is ultimately periodic.
\EThe
If $\a$ is an element of $\seqN$ or $\N^n$ ($n\in\N$), we will sometimes simply write $[\a]$ instead of $[a_1,\ldots]$ or $[a_1,\ldots,a_n]$ respectively.

\begin{figure}
 \includegraphics[scale=0.47]{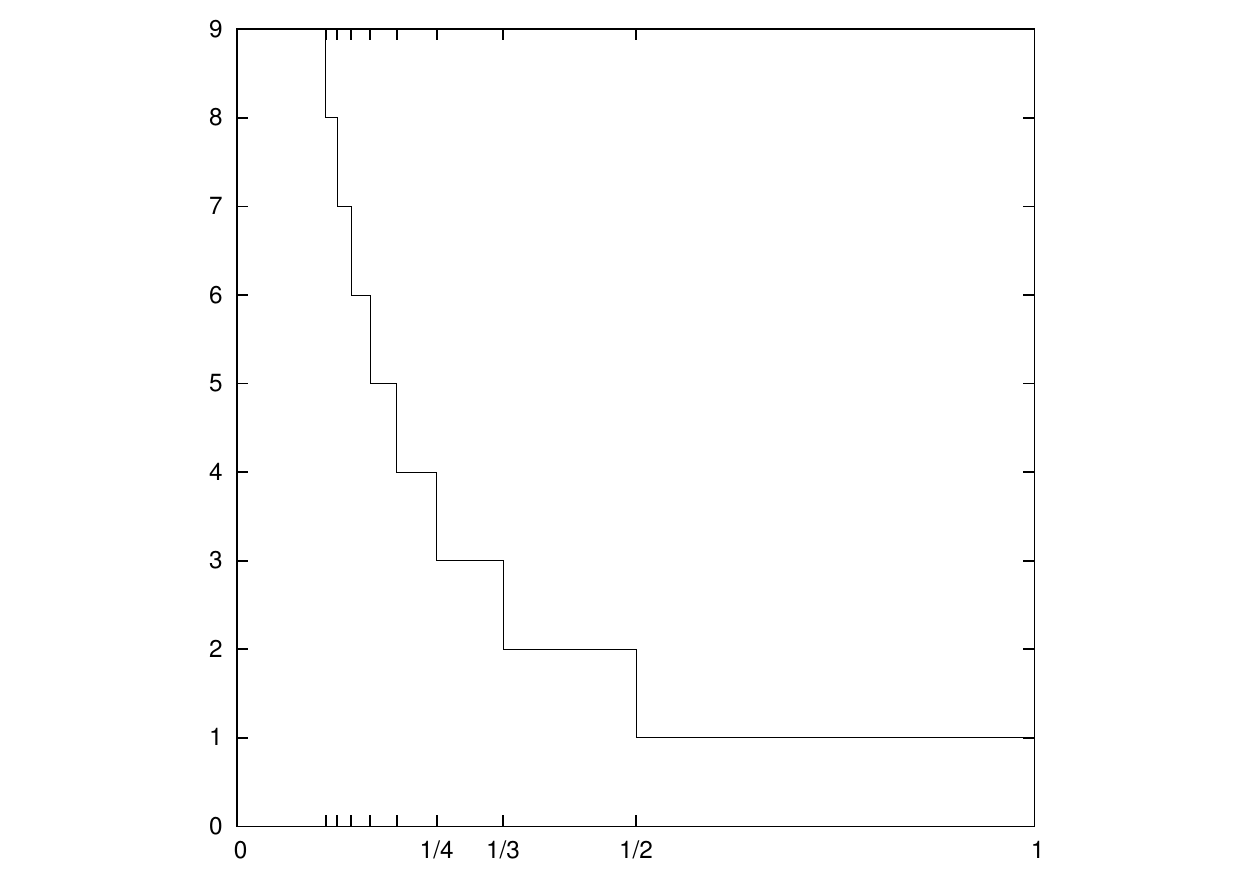} \includegraphics[scale=0.47]{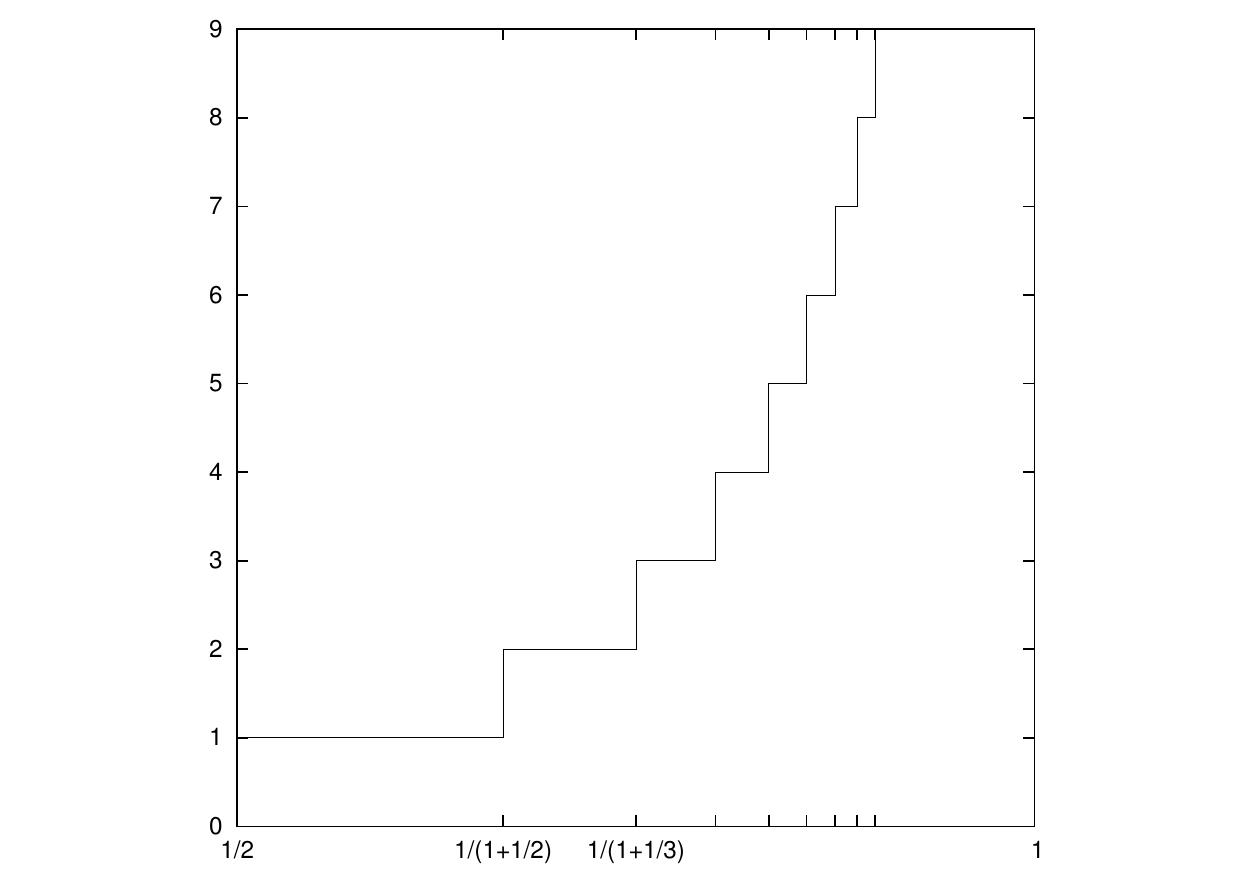}
\caption{The functions $x\mapsto a_1(x)$ (left panel) and $x\mapsto a_2(x)$ if $a_1(x)=1$ (right panel). This illustrates the fact that $I_1(x)$ is partitioned into a denumerably infinite number of intervals of rank $2$; in this case, $I_2(x)\subset [1/2,1]\cap I$, since $a_1(x)=1$ if and only if $x\in [1/2,1]\cap I$.}\label{fig:a_j}
\end{figure}

Let us now give a brief introduction to the notion of the metric theory of continued fractions. Since, for any
$\a\in\seqN$, $[\a]$ corresponds to an irrational number $x\in I$, one can consider, for each $j\in\N$, the term $a_j$ as a function of $x$ : $a_j=a_j(x)$. Let us fix $j\in\N$ and write $x=[a_1,\ldots,a_{j-1},r_j]$, with $r_j\in[1,\infty)$. It is easy to check that for any $k\in\N$, we have, if $j$ is odd, 
\[
 a_j=k
 \quad\text{if and only if}\quad
 \frac{1}{k+1} <r_j\le \frac{1}{k}
\]
and, if $j$ is even,
\[
 a_j=k
 \quad\text{if and only if}\quad
 k\le r_j < k+1.
\]
For any $j\in\N$, $a_j$ is thus a piecewise constant function. Moreover, $a_j$ is non-increasing if $j$ is odd and non-decreasing if $j$ is even. The functions $a_1$ and $a_2$ are represented in Figure~\ref{fig:a_j}. Let $x=[\a]$ be an irrational number; for $n\in\N$, we set
\[
 I_n(x)=\big\{y=[\b]\in I : b_j=a_j \text{ if } j\in\{1,\ldots,n\}\big\}.
\]
We will say that $I_n(x)$ is an interval of rank $n$. For any $n\in\N$, $I_n(x)$ is an irrational subinterval of $I$, $I_{n+1}(x)\subset I_n(x)$ and $\lim_n I_n(x)=\{x\}$. Indeed, using Proposition~\ref{pro:quotient avec reste} with $r_{n+1}=1$ and $r_{n+1}\to\infty$, one gets
\[
 I_n(x)= \left( \frac{p_n(\a)}{q_n(\a)}, \frac{p_n(\a)+p_{n-1}(\a)}{q_n(\a)+q_{n-1}(\a)} \right) \cap I,
\]
if $n$ is even (if $n$ is odd, the endpoints of the interval are reversed). Every interval of rank $n$ is partitioned into a denumerably infinite number of intervals of rank $n+1$. We will denote by $|I_n(x)|$ the Lebesgue measure of $I_n(x)$. One has, using Proposition~\ref{pro:convergents},
\[
 |I_n(x)| = \frac{1}{q_n(\a)(q_n(\a)+q_{n-1}(\a))}.
\]

\subsection{Cantor's bijection}
\begin{figure}
 \includegraphics[scale=0.47]{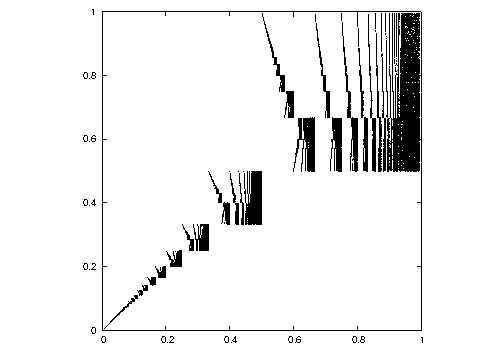} \includegraphics[scale=0.47]{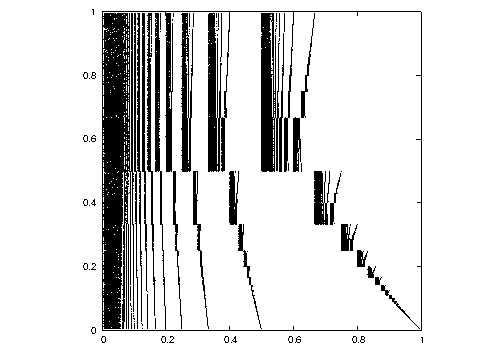}
\caption{The functions $f_1$ (left panel) and $f_2$ (right panel).}\label{fig:f}
\end{figure}

Cantor's bijection on $I$ (see \cite{cantor}) is a one-to-one mapping between $I$ and $I^2$. If $x\in I$, let $[a_1,\ldots]$ be the corresponding continued fraction and define the applications $f_1$ and $f_2$ as follows:
\[
 f_1(x)=[a_1,a_3,\ldots,a_{2j+1},\ldots]
 \quad\text{and}\quad
 f_2(x)=[a_2,a_4,\ldots,a_{2j},\ldots].
\]
These applications are represented in Figure~\ref{fig:f}.

Theorem~\ref{thm:expr cf} implies that the application
\[
 f: I\to I^2 \quad x\mapsto (f_1(x), f_2(x))
\]
is a one-to-one mapping. If $Q$ denotes the quadratic numbers of $I$, $f$ is a one-to-one mapping between $Q$ to $Q^2$. Since the cardinals of $E$ and $I$ are equal, $f$ can be extended to a one-to-one mapping from $E$ to $E^2$.

One can already show that Cantor's bijection is continuous. However, we will be more precise in the next section, using simpler arguments.
\BRem
For any $n\in\N$ and any $x\in I$, $f_1$ maps the interval $I_n(x)$ to $I_{m}(f_1(x))$, where $m=n/2$ if $n$ is even and $m=(n+1)/2$ if $n$ is odd. This indeed shows that $f_1$ is a continuous function; obviously, the same argument can be applied to $f_2$.
\ERem

\section{Continuity of Cantor's bijection on $I$}

Let $x\in I$; we write $\phi(x)=\a$ if $\a\in \seqN$ satisfies $x=[\a]$. For any $x\in \R$, $\floor{x}$ denotes the floor function and $\ceil{x}$ the ceil function: $\floor{x}=\sup\{ k\in\Z: k\le x \}$, $\ceil{x}= \inf \{k\in\Z : x\le k\}$.
\BPro\label{pro:phi homeo}
 The application $\phi$ is an homeomorphism between $I$ and $\seqN$.
\EPro
\BProof
 Let $x_j$ be a sequence on $I$ that converges to $x_0\in I$. The fact that $\phi(x_j)$ converges to $\phi(x_0)$ is a direct consequence of Euclid's algorithm, but it is even simpler when one has the metric theory of continued fractions at one's disposal. For any $n\in\N$, there exists $J\in\N$ such that $j\ge J$ implies $x_j\in I_n(x_0)$, which is sufficient.
%
%

Now let $\a_j$ be a sequence on $\seqN$ that converges to $\a_0\in \seqN$ and set $x_j=[a_{j,1},\ldots]$, $x_0=[a_{0,1},\ldots]$. For $\epsilon>0$, let $n\in\N$ such that
\[
 q_n(\a_0)(q_n(\a_0)+q_{n-1}(\a_0)) > \frac{1}{\epsilon}.
\]
Since there exists $J\in\N$ such that $x_j\in I_n(x_0)$ whenever $j\ge J$, one has
\[
|x_0-x_j|\le |I_n(x_0)|<\epsilon
\]
for such indexes. One can avoid the use of the metric theory of continued fractions using Proposition~\ref{pro:convergents}: for $\epsilon>0$, let $k\in\N$ such that $q_{k}(\a_0)>\sqrt{2/\epsilon}$. We have
\begin{eqnarray*}
 \lefteqn{|x_0-x_j|} &&\\
 &\le& |x_0- [a_{0,1},\ldots, a_{0,k}]| + |[a_{0,1},\ldots, a_{0,k}]-   [a_{j,1},\ldots, a_{j,k}]| \\
 && + | [a_{j,1},\ldots, a_{j,k}] -x_j| \\
 &\le& \frac{1}{q_k^2(\a_0)} + \frac{1}{q_k^2(\a_j)}
 = \frac{2}{q_k^2(\a_0)},
\end{eqnarray*}
for $j$ sufficiently large, which is sufficient to conclude.
\EProof
\BRem
We obviously have $[\cdot]=\phi^{-1}$ on $\seqN$.
\ERem
Since $(\seqN,d)$ is a separable complete metric space, we have reobtained the following well-known result.
\BCor
The space $I$ is a Polish space.
\ECor

\BPro
Cantor's bijection $f$ is an homeomorphism between $I$ and $I^2$.
\EPro
\BProof
 This is trivial since the application
\[
 \psi: \seqN \times \seqN \to \seqN
  \quad (\a,\b) \mapsto \c,
\]
where
\[
 c_j=\left\{\begin{tabular}{ll}
             $a_{(j+1)/2}$ & if $j$ is odd \\
             $b_{j/2}$      & if $j$ is even
            \end{tabular}\right.
\]
is an homeomorphism.
\EProof

Netto's theorem \cite{sagan} guarantees that such a function $f$ can not be extended to a continuous function from $E$ to $E^2$. The following result gives additional informations.
\BPro\label{pro:extension}
Any extension of Cantor's bijection to $E$ is discontinuous at any rational number.
\EPro
\BProof
 Let $x\in D$; there exists $k \in \N$ and $\a\in \N^k$ with $a_k>1$ such that
\[
 x=[a_1,\ldots,a_k]=[a_1,\ldots,a_k-1,1].
\]
Let $\b\in \seqN$ and set $x_j=[a_1,\ldots,a_k,r_j]$, $y_j=[a_1,\ldots,a_k -1,1,r_j]$ with $r_j=j+[\b]$. Both the sequences $x_j$ and $y_j$ converge to $x$ and it is easy to check that $\lim_j f(x_j)\not=\lim_j f(y_j)$.
\EProof

\section{H\"older regularity of Cantor's bijection}

In this section, we give some preliminary results about the H\"older regularity (see e.g.\ \cite{jaffard} and references therein) of Cantor's bijection.

Let $\alpha\in [0,1]$; a continuous and bounded real function $g$ defined on $A\subset\R$ belongs to the H\"older space $\Lambda^\alpha(x)$ with $x\in A$ if there exists a constant $C>0$ such that
\[
 |g(x)-g(y)|\le C |x-y|^\alpha,
\]
for any $y\in A$. The H\"older exponent $h_g(x)$ of $g$ at $x$ is defined as follows:
\[
 h_g(x)=\sup\{\alpha\in [0,1]: g\in \Lambda^\alpha(x)\}.
\]
If $h_g(x)<1$, then $g$ is not differentiable at $x$.

Let us now state our main result.
\BThe\label{thm:hol-reg}
Let $x=[\a]$ be an element of $I$ and $y\in I_n(x)\setminus I_{n+1}(x)$. One has
\[
 \frac{\displaystyle\frac{1}{n} \sum_{j=1}^{\ceil{n/2}} \log a_{2j-1}}{\displaystyle \frac{1}{n} \sum_{j=1}^{n+3} \log(a_j +1) +\frac{1}{n} C_1(n)}
 \le
 \frac{\log |f_1(x)-f_1(y)|}{\log |x-y|}
\]
and
\[
 \frac{\log |f_1(x)-f_1(y)|}{\log |x-y|}
 \le
 \frac{\displaystyle \frac{1}{n}\sum_{j=1}^{\ceil{n/2}+3} \log(a_{2j-1}+1)+\frac{1}{2n} C_2(n)}
{\displaystyle \frac{1}{n} \sum_{j=1}^n \log a_j},
\]
where
\[
 C_1(n)=\frac{\log 2}{2} + \log \max(\frac{a_{n+2}+2}{a_{n+2}+1},\frac{a_{n+3}+2}{a_{n+3}+1})
\]
and
\[
 C_2(n)=\frac{\log 2}{2} + \log \max(\frac{a_{2\ceil{n/2}+3}+2}{a_{2\ceil{n/2}+3}+1},\frac{a_{2\ceil{n/2}+5}+2}{a_{2\ceil{n/2}+5}+1}).
\]
\EThe
\BProof
Let $x=[\a]=[a_1,\ldots]$ be an element of $I$ and consider
\[
 y
 =[a_1,\ldots,a_n,b_{n+1},b_{n+2},\ldots],
\]
with $b_{n+1}\not=a_{n+1}$; for the sake of simplicity, one can suppose that $n$ is even. We will bound $|x-y|$ and $|f_1(x)-f_1(y)|$ with terms depending on $\a$ and $n$ only.

Since $I_n(x)=I_n(y)$, one has $|x-y|\le |I_n(x)|$ and
\[
 |I_n(x)| = \frac{1}{q_n^2(\a)} \frac{1}{1+q_{n-1}(\a)/q_n(\a)} \le \frac{1}{q_n^2(\a)}.
\]
Moreover, since
\begin{eqnarray*}
 q_n(\a)&=& a_n q_{n-1}(\a) + q_{n-2}(\a) \ge a_n q_{n-1}(\a) \\
 &\ge& a_n (a_{n-1} q_{n-2}(\a) + q_{n-3}(\a))
 \ge a_n\cdots a_3 (a_2 q_1(\a) +q_0(\a)) \\
 &\ge& a_n\cdots a_1,
\end{eqnarray*}
one gets
\[
 |x-y|\le \frac{1}{a_1^2 \cdots a_n^2}.
\]
The same reasoning can be applied to
\[
 f_1(x)=[a_1,a_3,\ldots,a_{n-1},a_{n+1},\ldots]
\]
and
\[
 f_1(y)=[a_1,a_3,\ldots,a_{n-1},b_{n+1},b_{n+3},\ldots]
\]
to obtain
\[
 |f_1(x)-f_1(y)|
 \le |I_{n/2}(f_1(x))|
 \le \frac{1}{a_1^2 a_3^2\cdots a_{n-1}^2}.
\]

For the lower bound of $|x-y|$, let us remark that $I_{n+1}(x)\cap I_{n+1}(y)=\emptyset$, but the distance between $I_{n+1}(x)$ and $I_{n+1}(y)$ can be zero. However, for any fixed $j\in\N$, there exists a denumerably infinite number of intervals of rank $n+1+j$ in between $I_{n+1+j}(x)$ and $I_{n+1+j}(y)$, i.e.\ there exists a denumerably infinite number of $z\in I$ such that $z'\in I_{n+1+j}(z)$ implies $x<z'<y$ or $y<z'<x$. If $z=[\c]$ is such an element, one has
\[
 |x-y|\ge |I_{n+3}(z)|
 \ge \frac{1}{q_{n+3}(\c)(q_{n+3}(\c)+q_{n+2}(\c))}
 \ge \frac{1}{2 q_{n+3}^2 (\c)}.
\]
The relations
\begin{eqnarray*}
 q_{n+3}(\c) &=& c_{n+3} q_{n+2}(\c) +q_{n+1}(\c)
 \le (c_{n+3}+1) q_{n+2}(\c) \\
 &\le& (c_{n+3}+1) (c_{n+2} q_{n+1}(\c) +q_{n}(\c))
 \le (c_{n+3}+1) \cdots (c_1+1)
\end{eqnarray*}
lead to
\[
 |I_{n+3}(z)| \ge \frac{1}{2(c_1+1)^2 \cdots (c_{n+3}+1)^2}. 
\]
Now let
\[
 j_0=\left\{\begin{tabular}{ll}
 $n+2$ & if $x<y$ \\
 $n+3$ & if $y<x$
 \end{tabular}\right.;
\]
one can choose $z$ such that $c_j=a_j$ for any $j\in\N$ except for the index $j_0$ for which $c_{j_0}=a_{j_0}+1$, so that $z>x$ in the case $x<y$  and $z<x$ in the case $y<x$. Moreover, $I_{n+1}(z)=I_{n+1}(x)\not=I_{n+1}(y)$, so that $x<z<y$ in the case $x<y$ and $y<z<x$ in the case $y<x$. One therefore has
\[
|x-y|\ge  |I_{n+3}(z)| \ge \frac{1}{2(a_1+1)^2 \cdots (a_{n+2}+1)^2 (a_{n+3}+2)^2},
\]
or
\[
|x-y|\ge  |I_{n+3}(z)| \ge \frac{1}{2(a_1+1)^2 \cdots (a_{n+2}+2)^2 (a_{n+3}+1)^2},
\]
depending on the value of $j_0$. Without loss of generality, one can assume that $j_0$ corresponds to the largest integer in such inequalities.

Now there also exists $w=[d_1,\ldots]$ such that $I_{n/2+3}(w)$ lies between $I_{n/2+3}(f_1(x))$ and $I_{n/2+3}(f_1(y))$; moreover one can choose $w$ such that $d_j=a_{2j-1}$ for any $j$ except for one index $j_0\in\{n/2+2, n/2+3\}$, for which $d_{j_0}=a_{2j_0 -1}+1$. One thus has
\[
 |f_1(x)-f_1(y)| \ge |I_{n/2+3}(w)| \ge \frac{1}{2(a_1+1)^2(a_3+1)^2\cdots (a_{n+3}+1)^2(a_{n+5}+2)^2}.
\]

Putting all these inequalities together and taking the logarithm, one gets
\[
 \frac{-2 \sum_{j=1}^{n/2} \log a_{2j-1}}{-\log 2 -2 \sum_{j=1}^{n+3} \log(a_j +1) -2\log(\frac{a_{n+3}+2}{a_{n+3}+1})}
 \le
 \frac{\log |f_1(x)-f_1(y)|}{\log |x-y|}
\]
and
\[
 \frac{\log |f_1(x)-f_1(y)|}{\log |x-y|}
 \le
 \frac{-\log 2- 2\sum_{j=1}^{n/2+3} \log(a_{2j-1}+1)-2\log(\frac{a_{n+5}+2}{a_{n+5}+1})}{-2 \sum_{j=1}^n \log a_j}.
\]
\EProof
Of course, the same reasoning can be applied to $f_2$, leading to the same result.
\BThe\label{thm:hol-reg2}
Let $x=[\a]$ be an element of $I$ and $y\in I_n(x)\setminus I_{n+1}(x)$. One has
\[
 \frac{\displaystyle \frac{1}{n} \sum_{j=1}^{\floor{n/2}} \log a_{2j}}{\displaystyle \frac{1}{n} \sum_{j=1}^{n+3} \log(a_j +1) +\frac{1}{n}C_1(n)}
 \le
 \frac{\log |f_2(x)-f_2(y)|}{\log |x-y|}
\]
and
\[
 \frac{\log |f_2(x)-f_2(y)|}{\log |x-y|}
 \le
 \frac{\displaystyle \frac{1}{n}\sum_{j=1}^{\floor{n/2}+3} \log(a_{2j}+1)+\frac{1}{n} C_2(n)}
{\displaystyle \frac{1}{n} \sum_{j=1}^n \log a_j},
\]
where $C_1$ is defined as in Theorem~\ref{thm:hol-reg} and
\[
 C_2(n)=\frac{\log 2}{2} + \log \max(\frac{a_{2\floor{n/2}+4}+2}{a_{2\floor{n/2}+4}+1},\frac{a_{2\floor{n/2}+6}+2}{a_{2\floor{n/2}+6}+1}).
\]
\EThe

To obtain a generic result about the regularity of Cantor's bijection, we need a direct consequence of the ergodic theorem on continued fractions \cite{ryll}. We say that a property $P$ concerning sequences of $\seqN$ holds almost everywhere if for almost every $x\in I$ (with respect to the Lebesgue measure), the sequence $\a\in \seqN$ such that $x=[\a]$ satisfies $P$. The following result can be obtained from the main theorem of \cite{nair}.
\BThe\label{thm:ergodic}
For any $k\in\N\cup\{0\}$, almost every sequence $\a\in\seqN$ satisfies
\[
 \frac{1}{n} \sum_{j=1}^n \log (a_j+k), \frac{1}{n} \sum_{j=1}^n \log (a_{2j}+k), \frac{1}{n} \sum_{j=1}^n \log (a_{2j-1}+k) \to \log K_k,
\]
as $n$ goes to infinity, where $K_k$ is defined by:
\[
 K_k=
\prod_{j=1}^\infty (1+ \frac{1}{j(j+2)})^{\log (j+k)/\log 2}.
\]
\EThe
The seminal result $\frac{1}{n} \sum_{j=1}^n \log a_j \to \log K_0$ was proven in \cite{khintchine} ; $K_0$ is called the Khintchine's constant. Here, we will be interested in the values
\[
 \log K_0\approx 0.987849056\cdots
 \qquad\text{and}\qquad
 \log K_1\approx 1.409785988\cdots
\]

Using Theorem~\ref{thm:hol-reg}, Theorem~\ref{thm:hol-reg2} and Theorem~\ref{thm:ergodic} as $n$ goes to infinity (or equivalently as $y$ tends to $x$), we get the following result.
\BCor\label{cor:fin cor}
For almost every $x\in I$, one has
\[
 h_{f_1}(x),h_{f_2}(x)\in \Big[\frac{\log K_0}{2\log K_1},\frac{\log K_1}{2\log K_0}\Big].
\]
\ECor
\BRem
 Let $\a\in \seqN$ be the sequence defined by
\[
 a_j=\left\{\begin{tabular}{ll}
  $2^j$ & if $j$ is even\\
  $1$   & if $j$ is odd
  \end{tabular}\right.,
\]
for any $j\in\N$ and set $x=[\a]$. It is easy to check, using Theorem~\ref{thm:hol-reg}, that for this particular point, we have $h_{f_1}(x)=0$, so that $f$ is a multifractal function.
\ERem

\BRem
The insiders of ergodic theory will certainly recognize the Birkhoff theorem (with the Gauss transformation, which preserves the Gauss measure and which is ergodic for this measure) behind some arguments to prove Corollary~\ref{cor:fin cor}.
\ERem

\section*{Acknowledgments}
The first author would like to thank S.~Jaffard for drawing his attention to this problem.


\begin{thebibliography}{10}

\bibitem{cantor}
G.~Cantor.
\newblock {E}in {B}eitrag zur {M}annigfaltigkeitslehre.
\newblock {\em Journal f{\"u}r die reine und angewandte Mathematik (Crelle's
  Journal)}, 1878(84):242--258, 1878.

\bibitem{gouvea}
F.Q. Gouv{\^e}a.
\newblock Was {C}antor surprised?
\newblock {\em Am. Math. Mon.}, 118(3):198--209, 2011.

\bibitem{jaffard}
S.~Jaffard.
\newblock Wavelet techniques in multifractal analysis.
\newblock In {\em Proceedings of Symposia in Pure Mathematics}, volume~72,
  pages 91--152, 2004.

\bibitem{jaffard-nicolay}
S.~Jaffard and S.~Nicolay.
\newblock Pointwise smoothness of space-filling functions.
\newblock {\em Applied and Computational Harmonic Analysis}, 26(2):181--199,
  2009.

\bibitem{jaffard-nicolay2}
S.~Jaffard and S.~Nicolay.
\newblock Space-filling functions and {D}avenport series.
\newblock {\em Recent Developments in Fractals and Related Fields}, pages
  19--34, 2010.

\bibitem{khintchine}
A.~Ya. Khintchine.
\newblock {\em Continued fractions}.
\newblock P. Noordhoff, 1963.

\bibitem{lothaire}
M.~Lothaire.
\newblock {\em Combinatorics on words}.
\newblock Cambridge University Press, 1997.

\bibitem{nair}
R.~Nair.
\newblock {O}n the {M}etrical {T}heory of {C}ontinued {F}ractions.
\newblock In {\em {P}roceedings of the {A}merican {M}athematical {S}ociety},
  volume 120, 1994.

\bibitem{ryll}
C.~Ryll-Nardzewski.
\newblock {O}n the {E}rgodic {T}heorems ({II}) : {E}rgodic {T}heory of
  {C}ontinued {F}ractions.
\newblock {\em {S}tudia {M}ath.}, 12:74--79, 1950.

\bibitem{sagan}
H.~Sagan.
\newblock {\em Space-filling curves}.
\newblock Springer, 1994.

\bibitem{wallace}
D.F. Wallace.
\newblock {\em Everything and more: a compact history of infinity}.
\newblock WW Norton \& Company, 2004.

\end{thebibliography}
\end{document}